\documentclass[12pt,amsfonts]{article}
\usepackage{graphicx}
\usepackage{latexsym}
\usepackage{amssymb}
\usepackage{amsmath}
\usepackage{enumerate}

\usepackage{layout}
\usepackage{eufrak}
\newtheorem{prop}{Proposition}
\newtheorem{lemma}{Lemma}
\newtheorem{definition}{Definition}
\newtheorem{corollary}{Corollary}
\newtheorem{theorem}{Theorem}
\newtheorem{remark}{Remark}

\def\real{{\mathord{{\rm I\kern-2.8pt R}}}}        
\def\inte{{\mathord{{\rm I\kern-2.8pt N}}}}
\def\div{{\mathrm{{\rm div}}}}

\def\sZZ{{\rm Z\kern-2.8ptem{}Z}}

\def\z{{\mathchoice
  {\sZZ}
  {\sZZ}
  {\rm Z\kern-0.30em{}Z}
  {\rm Z\kern-0.25em{}Z} }}
\def\sQQ{{\kern 0.27em \vrule height1.45ex width0.03em depth0em
          \kern-0.30em \rm Q}}
\def\qu{{\mathchoice
    {\sQQ}
    {\sQQ}
  {\kern 0.225em \vrule height1.05ex width0.025em depth0em \kern-0.25em \rm Q}
  {\kern 0.180em \vrule height0.78ex width0.020em depth0em \kern-0.20em \rm Q}
        }}
\def\sCC{{\kern 0.27em \vrule height1.45ex width0.03em depth0em
          \kern-0.30em \rm C}}
\def\complex{{\mathchoice
    {\sCC}
    {\sCC}
  {\kern 0.225em \vrule height1.05ex width0.025em depth0em \kern-0.25em \rm C}
  {\kern 0.180em \vrule height0.78ex width0.020em depth0em \kern-0.20em \rm C}
        }}

\newcommand{\R}{\mathbb{R}}
\newcommand{\N}{\mathbb{N}}

\newcommand{\ba}{\begin{array}}
\newcommand{\ea}{\end{array}}
\newcommand{\be}{\begin{equation}}
\newcommand{\ee}{\end{equation}}
\newcommand{\bea}{\begin{eqnarray}}
\newcommand{\eea}{\end{eqnarray}}
\newcommand{\beaa}{\begin{eqnarray*}}
\newcommand{\eeaa}{\end{eqnarray*}}

\def\z{\zeta}

\font\tenmath=msbm10 \font\sevenmath=msbm7 \font\fivemath=msbm5
\newfam\mathfam \textfont\mathfam=\tenmath
\scriptfont\mathfam=\sevenmath \scriptscriptfont\mathfam=\fivemath
\def\math{\fam\mathfam}

\def \={{\buildrel {\rm (law)} \over =}}
\def \N{{\math N}}
\def \R{{\math R}}

\def\qed{ \hfill \vrule width.25cm height.25cm depth0cm\smallskip}

\newcommand{\basa}{\begin{assumption}}
\newcommand{\easa}{\end{assumption}}

\newcommand{\bas}{\begin{assum}}
\newcommand{\eas}{\end{assum}}

\newcommand{\ignore}[1]{}
\begin{document}

\renewcommand{\thefootnote}{\fnsymbol{footnote}}

\renewcommand{\thefootnote}{\fnsymbol{footnote}}

\title{The density of the solution to the stochastic transport equation with fractional noise }

\author{Christian Olivera $^{1,}$\footnote{Supported by
FAPESP 2012/18739-0} \hskip0.2cm 
Ciprian A. Tudor $^{2,3,}$ \footnote{Supported by the CNCS grant PN-II-ID-PCCE-2011-2-0015. 
Associate member of the team Samos, Universit\'e de Panth\'eon-Sorbonne Paris 1 }\vspace*{0.1in} \\
$^{1}$ Departamento de Matem\'atica, Universidade Estadual de Campinas,\\
13.081-970-Campinas-SP-Brazil. \\
colivera@ime.unicamp.br \vspace*{0.1in} \\
 $^{2}$ Laboratoire Paul Painlev\'e, Universit\'e de Lille 1\\
 F-59655 Villeneuve d'Ascq, France.\\
$^{3}$Academy of Economical Studies, Bucharest, Romania \vspace*{0.1in}\\
 \quad tudor@math.univ-lille1.fr\vspace*{0.1in}}
\maketitle

\begin{abstract}
We consider the transport equation driven by the fractional Brownian motion. We study the existence and the uniqueness of the weak solution and, by using the tools of the Malliavin calculus, we  prove  the existence of the density of the solution and we give Gaussian estimates from above and from below for this density. 

\end{abstract}

 \medskip


 \smallskip


\section{Introduction}

The purpose of this paper is to study the probability law of the real-valued
solution of the following stochastic partial differential equations

\begin{equation}\label{trasportS}
 \left \{
\begin{array}{lll}
    du(t, x) + b(t, x) \nabla u(t, x) \ dt  +  \nabla u(t,x) \circ d B^{H}_{t} +  F(t,u) \ dt =0,\\
 u(0, x) = u_{0}(x),
\end{array}
\right .
\end{equation}
where $B^{H}_{t}=(B^{H_{1}}_{t},...,B^{H_{d}}_{t})$ is a fractional  Brownian motion (fBm) in
$\mathbb{R}^{d}$ with Hurst parameter $H=(H_{1},..., H_{d}) \in \left[ \frac{1}{2}, 1\right) ^{d}$  and the stochastic integration is understood in the symmetric (Stratonovich) sense. The equation (\ref{trasportS}) is usually called the stochastic transport equation and arises as a prototype model in  a wide variety of phenomena. Although we introduced (\ref{trasportS}) in a general form, we mention that some results will be obtained in dimension one.

The stochastic transport equation  with standard Brownian noise has been first studied in the celebrated works by Kunita \cite{Ku3}, \cite{Ku2} and  more recently it has been the object of study for many authors. We refer, among many others, to  \cite{CO}, \cite{Fre1}, \cite{FGP2}, \cite{Maurelli}, \cite{MNP}, \cite{Oli}.

Our aim is to analyze the stochastic partial equation (\ref{trasportS}) when the driving noise is the fractional Brownian motion, including the particular case of the Brownian motion. We will first give, by interpreting the stochastic integral in (\ref{trasportS}) as a symmetric integral via regularization in the Russo-Vallois sense \cite{RV},  an existence and uniqueness result for the weak solution to (\ref{trasportS}) via the so-called method of characteristics   and we express the solution as the initial value applied to the inverse flow generated by the equation of characteristics.  This holds, when $H_{i}=\frac{1}{2}, i=1,..,d$ for any dimension $d$ and in dimension $d=1$ if the Hurst parameter is bigger than one half. Using this representation of the solution to (\ref{trasportS}),  we study the existence and the Gaussian estimates for its density via the analysis of the dynamic of the inverse flow. A classical tools to study the absolute continuity of the law of random variables with respect to the Lebesque measure is the Malliavin calculus. We refer to the monographs \cite{N} or \cite{Sans} for various applications of the Malliavin calculus to the existence and smoothness of the density of random variables in general, and of solutions to stochastic  equations in particular.

 We will prove the Malliavin differentiability of the solution to (\ref{trasport}) by analyzing the dynamic of the inverse flow  generated by the characteristics (\ref{1}). Using a result in \cite{NV} we obtain, in dimension $d=1$ upper and lower Gaussian bounds for the density of the solution to the transport equation. We are also able to find the explicit form of the density in dimension $d\geq 2$ when the driving noise is the standard Brownian motion and the drift is divergence-free (i.e. the divergence of the drift vanishes).

We organized our paper as follows. In Section 2 we recall the existence and uniqueness results for the solution to the transport equation driven by the standard Brownian motion. In Section 3, we analyze the weak solution to the transport equation when the noise is the fBm, via the method of characteristics. In Section 4 we study the Malliavin differentiability of the solution to the equation of characteristics and this will be applied in Section 4 to obtain the existence and the Gaussian estimates for the solution to the transport equation. In Section 6 we obtain an explicit formula for the density when the noise is the Wiener process and the drift is divergence-free. 

\section{Stochastic transport equation driven by standard Brownian motion}

Throughout the paper, we will fix a probability space $(\Omega, \mathcal{F}, P)$ and a $d$-dimensional Wiener process $(B_{t}) _{t \in [0, T]}$ on this probability space. We will denote by $(\mathcal {F}_{t}) _{t\in [0,T]}$ the filtration generated by $B$. 

We will start by recalling some known facts on the solution to the transport equation driven by a standard Wiener process in $\mathbb{R} ^{d}$.

The equation (\ref{trasportS}) is interpreted in the strong sense, as the following stochastic
integral equation

\begin{equation}
u(t,x)=u_{0}(x)-\int_{0}^{t}b(s,x)\nabla u(s,x)\ ds -\sum_{i=0}^{d}\int_{0}%
^{t}   \partial_{x_i}  u(s,x)\circ dB_{s}^{i} \label{transintegral} -  \int_{0}^{t}   F(t,u) \ ds
\end{equation}
for $t\in [0,T]$ and $x\in \mathbb{R} ^{d}$.

The solution to (\ref{trasportS}) is related with the so-called equation of characteristics. That is, for $0\leq s\leq t$ and $x\in\mathbb{R}^{d}$, consider the following
stochastic differential equation in $\mathbb{R}^{d}$

\begin{equation}
\label{11}X_{s,t}(x)= x + \int_{s}^{t} b(r, X_{s,r}(x)) \ dr + B_{t}-B_{s}.
\end{equation}
and denote by $X_{t}(x): = X_{0,t}(x), t\in [0,T], x\in \mathbb{R} ^{d}$. 
 
For $m \in \N$ and $0< \alpha < 1$, let us assume the following hypothesis on $b$:
\begin{equation}
\label{REGULCLASS}
    b\in L^{1}((0,T); C_{b}^{m,\alpha}(\mathbb{R}^{d}))
\end{equation}
where $C^{m,\alpha}(\mathbb{R}^{d})$ denotes the class of functions of class $C ^{m}$  on $\mathbb{R}^{d}$ such that the last derivative is H\"older continuous of order $\alpha$.

Let us recall  the definition of the stochastic flow (see e.g. \cite{Ku}). 

\begin{definition}
 A stochastic flow is a family of maps $(\Phi _{s,t} : \mathbb{R} ^{d} \to \mathbb{R} ^{d} ) _{ 0\leq s\leq t\leq T}$ such that
\begin{itemize}

\item $\lim _{t\to s_{+} } \Phi_{s,t} (x) =x$ for every $x\in \mathbb{R} ^{d}$.

\item $\Phi _{u,t} \circ \Phi _{s,u} = \Phi _{s,t} $ if $0\leq s\leq u \leq t$.

\end{itemize} 
\end{definition}
Note that in \cite{Ku} the some measurability is also required in the definition of the flow but, since we are working  later in the paper with non-semimartingales, we will omit it.

\noindent It is well known that under conditions (\ref{REGULCLASS}), $X_{s,t}(x)$ is a
stochastic flow of $C^{m}$-diffeomorphism (see for example \cite{Ku2} and
\cite{Ku}). Moreover, the inverse $Y_{s,t}(x):=X_{s,t}^{-1}(x)$ satisfies the
following backward stochastic differential equation%

\begin{equation}
\label{itoassBac}Y_{s,t}(x)= x - \int_{s}^{t} b(r, Y_{r,t}(x)) \ dr - (B_{t}-B_{s}).
\end{equation}
for  every $0\leq s\leq t\leq T$, see \cite{FGP2} or \cite{Ku2} pp. 234.

In order to get the solution of (\ref{trasportS}) via the stochastic  characteristic method 
we considerer  the following ordinary differential equation 
\begin{equation}\label{deter}
Z_{t}(r)= r + \int_{0}^{t}   F(s,  Z_{s}(r)) \ ds .
\end{equation}

We have the following representation of the solution to the transport equation in terms of the inital data and of the inverse flow  (\ref{itoassBac}). We refer to e.g. \cite{Ku}  or \cite{Chow}, Section 3 for the proof.

\begin{lemma}
\label{lemaexis}  Assume (\ref{REGULCLASS}) for $m\geq3$ and let $u_{0}\in C^{m,\delta}(\mathbb{R}^{d}), F \in L^{\infty}((0,T); C ^{m}_{b}(\mathbb{R} ^{d})) $. Then the Cauchy problem (\ref{transintegral}) has a
unique solution $u(t,.)$ for $0\leq t\leq T$ such that it is a $C^{m}$-semimartingale which can be represented as%

\begin{equation*}
u(t,x)=Z_{t}(u_{0}(X_{t}^{-1}(x))),  \hskip0.3cm t \in [0,T], x\in \mathbb{R} ^{d}
\end{equation*}
where $Z$ is the unique solution to (\ref{deter}) and $X_{t} ^{-1} = X_{0,t} ^{-1}= Y_{0,t}$ for every $t\in [0,T]$.
\end{lemma}

\section{The weak solution of the transport equation driven by fractional Brownian motion}

We discuss in this section the existence, uniqueness and the representation of the solution to the standard equation driven by a fractional Brownian motion with Hurst parameter bigger than one half. We refer to the last section (the Appendix) for the basic properties of this process.   We will restrict throughout this section to the case $d=1$ and and we will use the concept of weak solution.  We explain at the end of this section (see Remark \ref{19a-2}) why we need to assume these restrictions.

 Consider the following  one-dimensional Cauchy problem: given an initial-data $u_0$,
find $u(t,x;\omega) \in \R$, satisfying 

\begin{equation}\label{trasport}
 \left \{
\begin{aligned}
    &\partial_t u(t, x;\omega) + \Big(\, \partial_x u(t, x;\omega)  \, \big(b(t, x) + \frac{d B ^{H}_{t}}{dt}(\omega)\big ) \Big)= 0,
    \\[5pt]
    &u|_{t=0}=  u_{0},
\end{aligned}
\right .
\end{equation}
with $T>0$, $\big( (t,x) \in U_T, \omega \in \Omega \big)$, where $U_T= [0,T] \times \R$, and $b:[0,T]+ \times
\R\to \R$ is a given vector field. The noise 
$B ^{H}$ is a fractional Brownian motion with Hurst parameter $H>\frac{1}{2}$ and  the stochastic
integral in (\ref{trasport}) will be understood  in the symmetric sense via regularization \cite{RV} or \cite{RV1}.  The fBm $B ^{H}$ is related to the Brownian moption $B$ via (\ref{BH}).

Let us first recall the notion of weak solution to (\ref{trasport}).

\begin{definition}\label{defisolu}  A stochastic process
$u\in  L^{\infty}(\Omega\times[0, T]\times \mathbb{R})$ is called 
a weak $L^{p}-$solution of the Cauchy problem \eqref{trasport},
when for any $\varphi \in C_c^{\infty}(\R)$,  $\int _{\mathbb{R}} u(t,
  x)\varphi(x)dx$ is an    adapted real value proces which has a continuous modification, finite covariation,  
 and for all $t \in [0,T]$, we have $P$-almost surely

\begin{eqnarray}
    \int_{\mathbb{R}} u(t,x) \varphi(x) dx &=& \int_{\mathbb{R}}  u_{0}(x) \varphi(x) \ dx
   +\int_{0}^{t} \!\! \int_{\mathbb{R}}  u(s,x) \ b(s,x) \partial_{x} \varphi(x) \ dx ds\nonumber 
\\
   &&\; + \int_{0}^{t} \!\! \int_{\mathbb{R}} u(s,x) \, b'(s,x) \, \varphi(x) \ dx ds \nonumber
\\
    &&\;  + \int_{0}^{t} \!\! \int_{\mathbb{R}} u(s,x) \ \partial_{x} \varphi(x) \ dx  d^{\circ}B ^{H} _s.\label{DISTINTSTR}
\end{eqnarray}
where $b'(s,x)$ denotes the derivative of $b(s,x)$ with respect to the variable $x$.
\end{definition}

\noindent At this point, we need to  recall the definition of the symmetric integral $d^{\circ} B ^{H}$ that appears in   (\ref{DISTINTSTR}). assume $(X_t)_{ t\geq 0} $ is  a continuous process and
$(Y_t)_{ t\geq 0}$ is  a process with paths in
$L_{loc}^{1}(\mathbb{R}^{+})$, i.e. for any $ b > 0$, $
\int_{0}^{b}|Y_t| dt <\infty$ a.s. The generalized stochastic
integrals (forward, backward and symmetric)  are defined through a regularization
procedure see \cite{RV}, \cite{RV1}. That is, let $ I^{0}(\epsilon, Y, dX)$  be the  $\varepsilon-$symmetric integral

$$
I^{0}(\epsilon, Y, dX)=\int_{0}^{t} Y_{s}
\frac{(X_{s+\epsilon}-X_{s-\epsilon})}{2\epsilon} ds \ t \geq 0.
$$
The symmetric integral $\int_{0}^{t} Y d^{\circ} X$ is defined as 

$$
   \int_{0}^{t} Y d^{\circ} X: =
\lim_{\epsilon\rightarrow 0}I^{0}(\varepsilon, Y, dX)(t),
$$
for every $t\in [0,T]$, provided the limit exist ucp (uniformly on compacts in probability).

Similarly to Lemma \ref{lemaexis}, we  also have a representation
formula for the weak solution  in terms of the initial condition $u_0$ and the (inverse) stochastic flow associated
to SDE  (\ref{1}).

\begin{theorem}\label{repre} Assume that  $b \in L^{\infty}((0,T); C_b^{1}(\mathbb{R}^{d}))$. Then there exists a $C^{1} (\mathbb{R}) $ stochastic flow of diffeomorhism $X_{s,t}, 0\leq s\leq t\leq T$ that satisfies 
\begin{equation}
\label{1}
X_{s,t} (x)= x+ \int_{s}^{t} b(u, X_{s,u} (x)) du + B ^{H} _{t}-B ^{H} _{s}
\end{equation}
for every $x\in \mathbb{R}^{d}$. 
Moreover, if $d=1$, given $u_{0}\in L^{\infty}(\mathbb{R})$, the stochastic process 
\begin{equation}
\label{19a-1}u(t, x):= u_0(X_{t}^{-1}(x)), \hskip0.5cm t\in [0,T], x\in \mathbb{R}
\end{equation} is the 
unique weak $L^{\infty}-$ solution of the Cauchy problem \eqref{trasport}, where
$X_{t}:= X_{0,t}$ for every $t\in [0,T]$.
\end{theorem}
{\bf Proof: } We will proceed in several steps: first we show that (\ref{1}) is a diffeomorphism flow, then we prove the uniquennes of the $L^{\infty}$ weak solution to (\ref{trasport}) and then we show that (\ref{19a-1}) satisfies the transport equation (\ref{trasport}). 

Let us first show that (\ref{1}) generates a flow of diffeomorphism. By doing the linear transformation 
$$Z_{s,t}=X_{s,t}(x)- ( B_{t}^{H}-B_{s}^{H} )$$ we deduce that the 
 equation (\ref{1}) is equivalent to the  random equation 

\begin{equation}\label{itoassequevalent}
Z_{s,t}(x)= x + \int_{s}^{t}   b(r, Z_{s,r}(x) + B_{r}^{H}-B_{s}^{H} ) \ dr  
\end{equation}for  $0 \leq s\leq t\leq T$.

From the classical theory for ordinary differential equations (see e.g. \cite{AMR}) we have that 
$Z_{s,t}(x)$ with $ 0\leq s\leq t\leq T$   is a $C^{1}(\mathbb{R}^{d}) $ diffeomorphism flow. Thus we deduce that 
$X_{s, t}(x)$   is a $C^{1}(\mathbb{R}^{d}) $ diffeomorphism flow.

In a second step, we will show that the transport equation with fBm noise admits a unique $L^{\infty}$ weak solution. By linearity we have to show that a weak
$L^{\infty}-$solution with initial condition $u_{0}=0$ vanishes
identically.  Applying the It\^o-Ventzel for the symmetric integral formula (see Proposition 9 of \cite{FlandRusso}) to  $F(y)=\int u(t,x) \varphi(x-y) \ dx $ (which depends on $\omega$), we obtain that

\begin{eqnarray}
\int _{\mathbb{R}} u(t,x) \varphi(x-B^{H}_{t}) dx
&=&
\int_{0}^{t} \int  _{\mathbb{R}} b(s,x) \partial_x  \varphi(x-B ^{H} _{s}) u(s,x)  dx
ds \nonumber\\
&&+ \int_{0}^{t} \int _{\mathbb{R}}  b'(s,x) \varphi(x-B ^{H} _{s}) u(s,x) 
dx ds \nonumber \\
&&
 +   \int_{0}^{t} \int _{\mathbb{R}} 
 u(s,x) \partial_x  \varphi(x-B ^{H}_{s}) dx d^{\circ} B ^{H}_{s}\nonumber \\
&& +   \int_{0}^{t} \int _{\mathbb{R}} 
 u(s,x) \partial_y  [ \varphi(x-B^{H}_{s})] dx  d^{\circ} B ^{H}_{s} \label{19a-4}.
\end{eqnarray}

\noindent We observe that $\partial_y [
\varphi(x-B^{H}_{s})]=- \partial_x \varphi(x-B ^{H}_{s})$. Thus the process 
$$V(t,x): =u(t,x+ {B ^{H}_{t}})$$ verifies

\begin{eqnarray*}
\int _{\mathbb{R}}V(t,x) \varphi(x) dx &=& \int_{0}^{t} \int _{\mathbb{R}}  b(s,x+B ^{H}_{s}) \partial_x 
\varphi(x) V(s,x)dx ds\\
&&
+ \int_{0}^{t} \int _{\mathbb{R}}  b'(s,x+B ^{H}_{s}) \varphi(x) V(s,x) dx ds.
\end{eqnarray*}

\noindent Let $\phi_{\varepsilon}$ be a standard mollifier and let 
$V_{\varepsilon}(t,x): =V(t,.)\ast \phi_{\varepsilon}$. Then it holds

 \begin{eqnarray*}
 \int _{\mathbb{R}} V(t,z) \phi_{\varepsilon}(x-z) dz&=&  
 \int_{0}^{t} \int _{\mathbb{R}} V(s,z) \, b(s,z+B ^{H}_{s}) \, \partial_z \phi_{\varepsilon}(x-z)   dz ds 
\\
&&+ \int_{0}^{t} \int _{\mathbb{R}}  u(s,z)  b'(s,z+B ^{H}_{s}) \, \phi_{\varepsilon}(x-z) dz ds 
\end{eqnarray*}

\noindent From an algebraic 
convenient manipulatio we get

\[
\frac{dV_{\varepsilon}}{dt}-b(t,x-B ^{H}_{t}) \partial_x  V_{\varepsilon}=\mathcal{R}_{\varepsilon}(b,u)
\]
where $ \mathcal{R}_{\varepsilon}(b,u)$ is the commutator defined as
$$
    \mathcal{R}_{\varepsilon}(b,u)=(b\partial _{x}) (\phi_{\varepsilon}\ast u )- \phi_{\varepsilon}\ast((b\partial _{x})u).
$$

 Since $b(s,x+B^{H}_{s})$ belongs   a.s to $ L^{\infty} ((0,T); C ^{1}_{b} (\mathbb{R} ))$  then by the Commuting Lemma (see Lemma II.1 of
\cite{DL}), the process  $V_{\varepsilon}(t,x)=V(t,.)\ast \phi_{\varepsilon}$
satisfies
\[
\lim_{\varepsilon \rightarrow
0}\frac{dV_{\varepsilon}}{dt}-b(t,x-B ^{H}_{t}) \partial_x  V_{\varepsilon}=0
 \mbox{ a.s.  in }\ L^{1}([0,T],
L^{1}_{loc}(\mathbb{R})).
\]

\noindent We deduce that if $\beta\in
C^{1}(\mathbb{R})$ and $\beta^{\prime}$ is bounded, then

\begin{equation}\label{norma}
\frac{d\beta(V)}{dt}-b(t,x-B ^{H}_{t}) \partial _{x} \beta(V)=0.
\end{equation}

\noindent Now, by  Theorem II.
2 of \cite{DL}, we define for each $M \in [0,\infty)$ the function
$\beta_M(t)=(|t| \wedge M )^p$ and obtain that
\[
\frac{d}{dt}\int \beta_M(V(t,x))dx \leq C \int \beta_M(V(t,x))dx.
\]
Taking expectation we have that
\[
\frac{d}{dt}\int \mathbb{E}(\beta_M(V(t,x)))dx \leq C \int
\mathbb{E}(\beta_M(V(t,x)))dx.
\]
From Gronwall Lemma we conclude that $\beta_M(V(t,x))=0$. Thus
 $u=0$.

Let us finally show that (\ref{19a-1}) satisfies (\ref{trasport}).  We have that, by the change of variables $X_{t} ^{-1}(y)=x$

\begin{equation}
\label{PUSHFORWARD}
    \int_{\mathbb{R}}  u_0(X_t^{-1}) (y) \; \varphi(y) dy
    =\int_{\R^d}  u_0(x) \;  X'_t (x)\varphi(X_t (x) )  dx ,
\end{equation}
for each $t \in [0,T]$, where $X' _{t}(x)$ denotes the derivative with respect to $x$ of $X_{t}(x)$.

Notice that $X'_{t}(x)= 1+ \int_{0} ^{t} b'(s, X_{s}(x)) X '_{s}(x)ds$ for every $t\in [0,T], x\in \mathbb{R}$. By applying  It\^o's formula (see \cite{RV1}, \cite{RV})  to  the product $$ X'_t (x) \varphi(X_t (x))$$  and using the fact that $B^{H}$ has zero quadratic variation when $H>\frac{1}{2}$
we obtain that

\begin{eqnarray}
  \int_{\mathbb{R}}  u_0 (X_t^{-1}(x))   \varphi(x)  dx  
     &=&  \int_{\mathbb{R}}   u_0(x)    dx  +\int_{0}^{t}  \int_{\mathbb{R}}   u_0(x)  b(s,X_s(x)) X'_t (x)  \cdot  \varphi ' (X_s (x))   dx  ds \\
&+& \int_{0}^{t}\int_{\mathbb{R}}  u_0(x)  b'(s,X_s (x) ) X'_t (x) \varphi '(X_s(x) )    dx  ds \nonumber \\
&+&  \int_{0}^{t} \int_{\mathbb{R}}  u_0(x)  X'_t   (x)  \varphi '(X_s(x)) dy   d^{\circ}B_{s}^{H}. \label{19a-3}
\end{eqnarray}
Note that the It\^o formula in \cite{RV1} guarantees the existence of the symmetric stochastic integrals in (\ref{19a-3}) above. Now, by the change variable $ y=X_t (x)$ we have that

\begin{eqnarray*}
 && \int_{\mathbb{R}}   u_0(X_t^{-1}(x))  \varphi(x)  dx  
     =  \int_{\mathbb{R}}   u_0(x)    dx  + \int_{0}^{t} \int_{\mathbb{R}}  u_0(X_s^{-1}(x))  b(s,y)   \cdot  \varphi '(y)   dy  ds \\
&&
+\int_{0}^{t} \int_{\mathbb{R}}  u_0(X_s^{-1}(x))    b'(s,y)    \varphi'(y)    dy  ds 
\\
&&+ \int_{0}^{t} \int_{\mathbb{R}}   u_0(X_s^{-1}(x))  \;    \varphi '(y) dy    d^{\circ} B_{s}^{H}.
\end{eqnarray*}

From this we conclude que $u(t,x)=u_0(X_t^{-1}(x))$ is a weak solution of (\ref{trasport}). Its adaptedness is a consequence of (\ref{BH}). Thus the  unique solution  to (\ref{trasport}) is $u(t,x)=u_0(X_t^{-1}(x))$ for every $t\in [0,T]$ and for every $x\in \mathbb{R}$.  \qed

\begin{remark}
\label{19a-2}
\begin{itemize}
\item We need to restrict to the situation $d=1$ in order to get the existence of the symmetric integral in (\ref{19a-3}) or (\ref{19a-4}). Here we also used the hypothesis $H> \frac{1}{2}$ that ensures that there is not a second derivative term in the It\^o formula.

\item The uniqueness of the weak solution can be obtined with weaker assumption on the drift  $b$ by following the proof of Theorem 3.1 in \cite{CO}.
\end{itemize}

\end{remark}

\section{Fractional Brownian flow}

In this section we will analyze the properties of the stochastic flow generated by the fractional Brownian motion. We will call it the fractional Brownian flow in the sequel.  
Fix $d\geq 1$ and let $ B^{H}= ( B ^{H_{1}}, B^{H_{2}}, \ldots, B^{H_{d}})$ be a $d$-dimensional fractional Brownian motion with Hurst parameter $H= (H_{1}, H_{2}, \ldots , H_{d}) \in (0,1) ^{d}$.

Recall (see Theorem \ref{repre}) that if $b \in L^{\infty}((0,T), C^{1}_{b} (\mathbb{R} ^{d}))$, (\ref{1}) generates a $C^{1}$-stochastic flow of diffeomorphism. We next describe the dynamic of the inverse flow of (\ref{1}).

\begin{lemma}
Let $b \in L^{\infty}((0,T), C^{1}_{b} (\mathbb{R} ^{d}))$ and denote, for every $0\leq s\leq t\leq T $ and for every $x\in \mathbb{R} ^{d}$
\begin{equation}
\label{y}
Y_{s,t} (x)= X_{s, t} ^{-1}(x)
\end{equation}
the inverse of the stochastic flow given by (\ref{1}). Then the inverse flow satisfies the backward stochastic equation
\begin{equation}
\label{back}
Y_{s,t}(x)= x-\int_{s}^{t} b(r, Y_{r,t} )dr -( B^{H}_{t} -B^{H}_{s})
\end{equation}
for every $x\in \mathbb{R} ^{d}.$
\end{lemma}
{\bf Proof: } It follows from Kunita \cite{Ku}. Indeed, Lemms 6.2, page 235 in \cite{Ku} says that for any continuous function in two variables $g$ we have
\begin{equation*}
\int_{s}^{t} g(r, X_{s,r}(y)) dr | _{y=X_{s,t}^{-1} (x) } = \int_{s}^{t} g(r, X_{r,t}^{-1} (x))dr
\end{equation*}
and it suffices to apply the above identity to (\ref{1}). \qed

We need the following useful lemma.
\begin{lemma} \label{6a-5}
Let us introduce the notation, for $t\in [0,T]$ and $x\in \mathbb{R} ^{d}$,
\begin{equation}
\label{r}
R_{t,x} (u)= Y_{ t-u, t} (x), \hskip0.5cm \mbox{ if } u \in [0,t].
\end{equation}
Then we have, for every $t\in [0,T], u\in [0,t] $ and $x\in \mathbb{R} ^{d}$
\begin{equation}
\label{5a-1}
R_{t,x}(u)= x- \int_{0} ^{u} b(a, R_{t, x} (a) )da -( B ^{H}_{t}- B ^{H} _{t-u}).
\end{equation}
\end{lemma}
{\bf Proof: } In (\ref{y}) we use the change of notation $u=t-s$ and we get for every $y\in \mathbb{R} ^{d}$,
$$ R_{t,y}(u)= y- \int_{t-u}^{t} b( r, Y _{r,t} (x))dr -(B_{t}^{H}- B_{t-u}^{H})$$
and then,  with the change of variables $a=t-r$ in the integral $dr$, we can writye
$$R_{t,y}(u)= y-\int_{0} ^{u} b(a, Y_{t,y}(a) )da -( B ^{H}_{t}- B ^{H}_{t-u})$$
with $R_{u,y}(u)=y$. \qed

As a consequence of the above Lemma \ref{6a-5}, we get the uniqueness of solution to the backward equation (\ref{back}) satisfied by the inverse flow.

\begin{corollary}
If $ (\tilde{Y} _{s,t} ) _{0\leq s\leq t\leq T}$ is another two parameter process that satisfies (\ref{y}) with $\tilde{Y} _{s,s}(x)=x$ and $b $ is Lipschitz in $x$ uniformy with respect to $t$, then $\tilde{Y}_{s,t} (x)= Y_{s,t}(x) $ for every $0\leq s\leq t$ and for every $x\in \mathbb{R} ^{d}$.
\end{corollary}
{\bf Proof: } If $\tilde{Y}$ satisfies (\ref{y}), then, if we denote $\tilde{R}_{t,x}(u)=\tilde{ Y}_{t-u,t}(x) $, we get from Lemma that $\tilde{R}$ satisfies (\ref{r}) and the Gronwall lemma and the Lipschitz assumption on the drift $b$ imply the conclusion. \qed

We denote by $D$ the Malliavin derivative with respect with the fBm $ B^{H}$ (see the Appendix).  
\begin{prop}\label{5a-3}
Assume $b \in L^{\infty}((0,T), C^{1}_{b} (\mathbb{R} ^{d}))$ and let $X_{s,t}$ be given by (\ref{1}). Then, for every $0\leq s\leq t \leq T$ and for every $x\in \mathbb{R}^{d}$,  the components of inverse flow $Y^{i}_{s,t} $ ($1\leq i\leq d$) are  Malliavin differentiable and for every $\alpha \in [s, t]$
\begin{equation*}
D_{\alpha}  Y^{i}_{s,t}(x)=-\int_{s}^{t} \sum_{j=1} ^{d} \frac{\partial b^{i}}{\partial x_{j}}(r, Y_{r,t}) D_{\alpha} Y ^{j} _{r,t} (x) dr -1
\end{equation*}
and $D_{\alpha}  Y^{i}_{s,t}(x)=0$ if $\alpha \notin [s,t]$. We denoted by $b^{i}$ ($1\leq i\leq d$) the components of the vector mapping $b$.
\end{prop}
{\bf Proof: } It suffices to show that the random variable $R_{t,x}(u)$ defined by (\ref{r}) is Malliavin differentiable for any $x\in \mathbb{R} ^{d}$ and for every $0\leq u\leq t\leq T$.  We will give the sketsch of the proof which  follows by a  routine fix point argument. Fix $x\in \mathbb{R} ^{d}, t \in [0,T]$ and define the iterations
$$R^{(0)} _{t,x}(u)=x, \mbox{ for every } u\in [0,t]$$
and for $n\geq 1$,
\begin{equation*}
R_{t,x} ^{(n)}(u)= x-\int_{0} ^{u} b(a, R^{(n-1)}_{t, x} (a) )da -( B ^{H}_{t}- B ^{H} _{t-u}).
\end{equation*}
By induction, we can prove by standard arguments (see e.g. \cite{N}, Theorem 2.2.1) that for every $p\geq 1$
$$\sup_{0\leq u\leq t} \mathbf{E} \left| R ^{(n) }_{t,x} (u)\right| ^{p} <\infty, $$

$$ R^{(n), j} _{t,x}(u) \in \mathbb{D} ^{1, \infty}, \hskip0.5cm j=1,.., d, $$
and 
$$\sup _{n \geq 1}\sup_{ \alpha \in [0,T]} \mathbf{E} \left| D_{\alpha}  R^{(n), j} _{t,x}(u) \right| ^{p} < \infty$$
where 
$ R^{(n), j} _{t,x}(u) $ denotes the $j$ th component of $ R^{(n)} _{t,x}(u) .$ Moreover, the sequence of random variables  $ (R^{n} _{t,x}(u)) _{n \geq 1}$  converges in $L^{p}$ to $R_{t,x}(u)$ which is the unique solution to (\ref{5a-1}). 
It follows from Lemma  1.2.3 in \cite{N} that $R _{t,x}(u)$ belongs to $ \mathbb{D} ^{1,\infty}$.
\qed

\begin{remark}
Note that, when the noise is the standard Brownian motion, the Malliavin differentiability of $Y$ is also claimed in \cite{MNP}.
\end{remark}

\section{Existence and Gaussian bounds for  the  density of the solution to the transport equation in dimension one} 
In this section we will assume that $d=1$.  On the other hand, the results in these section (except Theorem \ref{6a-3}) will hold for every $H\in (0,1)$. We also mention that we will use the notation $c,C..$ for generic positive constants that may vary from line to line.

From Proposition \ref{5a-3} we immediately obtain the explicit expression for the Malliavin derivative of the inverse flow.

\begin{prop}
If $b \in L^{\infty}((0,T), C^{1}_{b} (\mathbb{R} ^{d}))$ and $Y_{s,t}$ is defined by (\ref{y}), we have for every $\alpha$ and for every $0\leq s\leq t\leq T$
\begin{equation}
\label{5a-4}
D_{\alpha} Y _{s,t} (x)=- 1_{[s,t]} (\alpha )e ^{- \int_{s} ^{\alpha} b'(r, Y_{r,t}) dr}
\end{equation}
with $b'(t,x)$ the derivative of $b(t,x)$ with respect to $x$.
\end{prop}
{\bf Proof: } For every $\alpha $, we have
\begin{eqnarray*}
D_{\alpha} Y_{s,t} (x) &=& -\int_{s}^{t} b'(s_{1}, Y_{s_{1}, t} (x)) D_{\alpha } Y _{s_{1}, t}(x) ds_{1}- 1_{[s,t] } (\alpha) 
\end{eqnarray*}
and by iterating the above relation we can write, for every $0\leq s\leq t\leq T$ and for evey $\alpha \in [0,T]$,
\begin{eqnarray*}
D_{\alpha} Y_{s,t} (x) &=& -1_{[s,t] } (\alpha ) \sum_{n\geq 0} (-1) ^{n} \int_{s}^{\alpha} ds_{1} \int_{s_{1} }^{\alpha} ds_{2} ...\int_{s_{n-1} }^{\alpha} ds_{n}\\
&&\times  b'(s_{1}, Y_{s_{1}, t} (x))b'(s_{2}, Y_{s_{2}, t} (x))..b'(s_{n}, Y_{s_{n}, t} (x))\\
&=& -1_{[s,t] } (\alpha ) \sum_{n\geq 0}\frac{ (-1) ^{n} }{n!} \left( \int_{s} ^{t} dr b'(r, Y_{r,t}(x) )\right) ^{n} \\
&=& - 1_{[s,t]} (\alpha )e ^{- \int_{s} ^{\alpha} b'(r, Y_{r,t}) dr}.
\end{eqnarray*}
\qed

The main tool in order to obtain the Gaussian estimates for the density of the solution to the trasport equation is the following result given in \cite{NV}.

\begin{prop}\label{6a-4}

I If $ F \in \mathbb{F} ^{1,2}$, let 
$$g_{F}(F)= \int_{0} ^{\infty} d\theta e^{-\theta} \mathbf{E}  \left[ \mathbf{E}  '\left( \langle DF, \widetilde{DF} \rangle _{\mathcal{H}} | F \right) \right]$$
where for any random variable $X$, we denoted 
$$\tilde{X}(\omega , \omega ') = X ( e ^{-\theta }w + \sqrt{1-e^{-2\theta} }\omega ').$$
Here  $\tilde{X}$ is  defined on a product probability space
 $\left( \Omega \times \Omega ', {\cal{F}} \otimes {\cal{F}}, P\times P'\right)$ and ${\mathbf E}'$
 denotes the expectation with respect to the probability measure $P'$. 
If there exists two constants $\gamma _{min}$ and $\gamma _{max}$ such that  almost surely 
$$0\leq \gamma_{ min} \leq g_{F}(F) \leq \gamma _{max} $$ 
then $F$ admits a density $\rho$. Moreover, for every $z\in \mathbb{R}$,
\begin{equation*}
\frac{ \mathbf{E} \vert F -\mathbf{E} F \vert }{2\gamma ^{2}_{max} }e ^{ -\frac{(z-\mathbf{E} F) ^{2}}{2\gamma^{2} _{min} } }\leq \rho(z) \leq 
\frac{ \mathbf{E} \vert F -\mathbf{E} F \vert }{2\gamma^{2} _{min} }e ^{ -\frac{(z-\mathbf{E} F) ^{2}}{2\gamma ^{2} _{max} }}
\end{equation*}

\end{prop}

To apply the above result, we need to controll the Malliavin derivative of the inverse flow. This will be done in the next result. Notice that a similar method has been used in e.g. \cite{AB}, \cite{BKT} or \cite{NQ} for various types of stochastic equations. In the sequel $\mathcal{H}$ denotes the canonical Hilbert space associated to the fractional Brownian motion (see the Appendix). 

\begin{prop}\label{6a-3}
Assume $H>\frac{1}{2}$ and $b \in L^{\infty}((0,T); C^{1}_{b}(\mathbb{R}))$.  Then there exist two positive constants $c<C$ such that for every $t\in [0,T]$ and for every $x \in \mathbb{R } $
\begin{equation}\label{6a-2}
ct^{2H} \leq \langle  D Y_{0,t}(x), \widetilde{ D Y_{0,t}(x)} \rangle  _{\mathcal{ H}}  \leq Ct ^{2H }
\end{equation}
where $Y_{0,t}$ is given by (\ref{back}).

\end{prop}
{\bf Proof: } Assume $H=\frac{1}{2}$. Then $H= L ^{2} ([0,T]) $ and 
$$ \langle  D Y_{0,t}(x), \widetilde{ D Y_{0,t}(x)} \rangle  _{\mathcal{ H}}= \int_{0} ^{t} d\alpha  e ^{- \int_{0} ^{\alpha} b'(r, Y_{r,t}(x)) dr}  e ^{- \int_{0} ^{\alpha} b'(r, \widetilde{Y_{r,t}(x)}) dr}  $$
and since
\begin{equation}
\label{6a-1}
e^{-T\Vert b'\Vert _{\infty} } \leq e ^{- \int_{s} ^{\alpha} b'(r, Y_{r,t}) dr}\leq e ^{T\Vert b'\Vert _{\infty}}
\end{equation} (and a similar bound holds for the tilde process) we obtain
$$ct \leq \langle  D Y_{0,t}(x), \widetilde{ D Y_{0,t}(x)} \rangle  _{\mathcal{ H}}   \leq Ct $$
with two positive constant $c$ and $C$.

Assume $H> \frac{1}{2}$. Then  by (\ref{27i-1})
$$ \langle  D Y_{0,t}(x), \widetilde{ D Y_{0,t}(x)} \rangle  _{\mathcal{ H}}   = \alpha _{H}\int_{0} ^{t} d\alpha \int_{0} ^{t} d\beta e ^{- \int_{0} ^{\alpha} b'(r, Y_{r,t}(x)) dr}e ^{- \int_{0} ^{\beta } b'(r, \widetilde{Y_{r,t}(x)}) dr}\vert \alpha - \beta \vert ^{2H-2} $$
and inequality (\ref{6a-1}) implies that
$$c\int_{0} ^{t} d\alpha \int_{0} ^{t} d\beta \vert \alpha - \beta \vert ^{2H-2} \leq \langle  D Y_{0,t}(x), \widetilde{ D Y_{0,t}(x)} \rangle  _{\mathcal{ H}}   \leq C\int_{0} ^{t} d\alpha \int_{0} ^{t} d\beta \vert \alpha - \beta \vert ^{2H-2} $$
which immediately gives (\ref{6a-2}).

Assume $H< \frac{1}{2}$. Then Proposition \ref{6a-1} in \cite{BKT} implies the lower bound in (\ref{6a-2}). Concerning the upper bound, it suffices again to follow \cite{BKT}, Section 3.4 and to note that for every $\alpha, \beta \in [0,T]$ with $\alpha >\beta$ we have
$$\vert D_{\alpha } Y _{0, t} (x) -D_{\beta} Y _{0,t} (x) \vert \leq e ^{- \int_{0} ^{\beta } b'(r, Y_{r,t}) dr}\left| e ^{- \int_{\beta } ^{\alpha} b'(r, Y_{r,t}) dr}-1 \right| \leq c(\alpha -\beta)$$
and the same bound holds for the $\widetilde{Y_{0,t}}$.
\qed

\vskip0.3cm

Denote by $m: =\mathbf{E}  u(t,x)$ (it satisfies a parabolic equation, see e.g. \cite{FGP2}).  We are ready to state our main result.

\begin{theorem}
Let $u(t,x)$ be the solution to the transport equation (\ref{trasport}). Assume that $u_{0} \in C^{1}(\mathbb{R})$ such that there exist $0<c<C$ with $c\leq u_{0}'(x) <C$ for every $x\in \mathbb{R}$ and $b\in L^{\infty}((0,T); C^{1}_{b}(\mathbb{R}))$. 
Then, for every $t\in [0,T]$ and for every $x\in \mathbb{R}$, the random variable $u(t,x)$ is Malliavin differentiable. Moreover $u(t,x)$ admits a density $\rho _{u(t,x) }$ and  there exist two positive constants $c_{1}, c_{2}$ such that 
\begin{equation}
\label{5a-5}
\frac{\mathbf{E} \vert u(t,x)- m\vert }{2c_{1} t^{2H}}e ^{-\frac{(y-m) ^{2}}{2c_{2} t^{2H}}}\leq \rho_{ u(t,x)} \leq  \frac{\mathbf{E} \vert u(t,x)- m\vert }{2c_{2} t^{2H}}e ^{-\frac{(y-m) ^{2}}{2c_{1} t^{2H}}}
\end{equation}

\end{theorem}
{\bf Proof: } Since by Theorem \ref{repre}, $u(t,x)= u_{0} (Y_{0,t}(x))$, we get the Malliavin differentiability of $u(t,x)$ from Proposition \ref{5a-3} and the chain rule for the Malliavin derivative  (see e.g. \cite{N}). Moreover, the chain rule  implies
$$D_{\alpha} u(t,x)= u_{0} '( Y_{0,t} (x))  D_{\alpha } Y_{0,t} (x) $$
and thus
$$ \langle  D u(t,x) , \widetilde{Du(t,x)} \rangle  _{ \mathcal{H}}=  u_{0} '( Y_{0,t} (x)) \widetilde{ u_{0} '( Y_{0,t} (x))} \langle  D Y_{0,t}(x), \widetilde{ D Y_{0,t}(x)} \rangle  _{\mathcal{ H}}  .$$
By Proposition \ref{6a-3} and the asumption $u_{0} \in C^{1}_{b}$, there exists two strictly positive constants $c<C$ such that 
$$ct^{2H} \leq  \langle  D u(t,x) , \widetilde{Du(t,x)} \rangle  _{ \mathcal{H}}\leq C t^{2H} $$
for every $t\in [0,T]$ and for every $x \in \mathbb{R}$. Now, Proposition \ref{6a-4}, point 2. implies that, if $F=u(t,x)$ then
$$c_{1} t ^{2H} \leq g_{F}(F) \leq c_{2} t^{2H}$$
and Proposition \ref{6a-4}, point 1. gives the conclusion. \qed

\section{Explicit expression of the density when the noise is the Brownian motion in $\mathbb{R} ^{d}$}

We obtained above the existence and Gaussian estimate for the solution to the transport equation in dimension 1 and for $H\geq \frac{1}{2}$. 
In this section, we will assume $d\geq 2$, $H=\frac{1}{2}$, that is, the transport equation is driven a standard Brownian motion in $\mathbb{R} ^{d}$. We obtain the followin explicit expression for the density of the solution when the divergence of the drift $b$ vanishes.

\begin{theorem}\label{um} Assume $d\geq 2$and let $u_{0}$ be a $ C^{m,\delta}(\mathbb{R}^{d})  $ diffeomorphism.  Assume 
(\ref{REGULCLASS}) for $m\geq3$. Moreover, suppose that
\begin{equation}
\label{DIVB}
  div \,  b= 0 . 
\end{equation}
 Fix $t \in [0,T]$  and  $x\in\mathbb{R}^{d} $. Then  the law 
of the   solution of   (\ref{trasportS}),  has a density  $\widetilde{\rho}$ with respect to the Lebesgue measure. Moreover
the density $\widetilde{\rho}$  admits the representation

\begin{equation}\label{9a-2}
 \widetilde{\rho}= Ju_{0}( Z_t^{-1}(y)) JZ_{t}  \rho(u_{0}^{-1}( Z_t^{-1}(y) ),t,x) 
\end{equation}
where $ \rho$ denotes the density of  the  solution to  (\ref{11}). 
\end{theorem}
{\bf Proof: } Let $u(t,x)$ solution of the SPDE    (\ref{trasportS} ). By  Lemma \ref{lemaexis}
we have that $u(t,x)$ has the representation

\begin{equation*}
u(t,x)=Z_t(u_{0}(X_{t}^{-1}(x))). 
\end{equation*}

 Let $\phi_{\varepsilon}$ be a standard mollifier and consider a smooth function $\varphi \in C_c^{\infty}( \mathbb{R}^{d} )$. Then

\begin{eqnarray*}
\mathbf{E} [\varphi(u(t,x))]&=& 
\mathbf{E} [\varphi(Z_t(u_{0}(X_{t}^{-1}(x))) )]\\
&=&
\mathbf{E}  \ [  \lim_{\epsilon\rightarrow 0 }  \int_{\mathbb{R}^{d}}   \phi_{\varepsilon}(y-x)    \varphi(Z_t(u_{0}(X_{t}^{-1}(y))) )   dy ] \\
&=&
 \lim_{\epsilon\rightarrow 0 }   \mathbf{E}\  [  \int_{\mathbb{R}^{d}}   \phi_{\varepsilon}(y-x)    \varphi(Z_t(u_{0}(X_{t}^{-1}(y))) )   dy ].
\end{eqnarray*}

The assumption  (\ref{DIVB}) implies that  $JX_t=1$, where $JX_t$ denote of the Jacobian map of $X_t$.  By doing one more time a chamge of variable, we can write

\begin{eqnarray}
\mathbf{E} [\varphi(u(t,x))]&=& 
 \lim_{\epsilon\rightarrow 0 }   \mathbf{E} [  \int_{\mathbb{R}^{d}}   \phi_{\varepsilon}(y-x)    \varphi(Z_t(u_{0}(X_{t}^{-1}(y))) )   dy ]\nonumber \\
&=&
 \lim_{\epsilon\rightarrow 0 }   \mathbf{E} [  \int_{\mathbb{R}^{d}}   \phi_{\varepsilon}(X_t(y)-x)    \varphi(Z_t(u_{0}(y)) )   dy ] \nonumber \\
&=&
 \lim_{\epsilon\rightarrow 0 }   \int_{\mathbb{R}^{d}}    \mathbf{E} [ \phi_{\varepsilon}(X_t(y)-x) ]   \varphi(Z_t(u_{0}(y)) )   dy ] \label{9a-1}
\end{eqnarray}
The random variable $X_{t}(x) $ admits a density $\rho$ in any dimension $d$. This is an easy consequence of equation (\ref{1}) (see e.g. \cite{N}). Therefore, (\ref{9a-1}) becomes
\begin{eqnarray*}
\mathbf{E} [\varphi(U(t,x))]&=&
 \lim_{\epsilon\rightarrow 0 }   \int_{\mathbb{R}^{d}}    \mathbf{E} [ \phi_{\varepsilon}(X_t(y)-x) ]   \varphi(Z_t(u_{0}(y)) )   dy ]\\
&=&
 \lim_{\epsilon\rightarrow 0 }   \int_{\mathbb{R}^{d}}    \int_{\mathbb{R}^{d}} \phi_{\varepsilon}(u-x)  \rho(u,t,y) du  \ \varphi(Z_t(u_{0}(y)) )   dy 
\end{eqnarray*}
and by calculating the limit above when $\varepsilon \to 0$ we get
\[
\mathbf{E} [\varphi(u(t,x))]=   
 \int_{\mathbb{R}^{d}}   \rho(y,t,x) \varphi(Z_t(u_{0}(y)) )   dy .
\]
Finally, the making succesively the changes of variables $w=u_{0} (y)$ and  $y= Z_{t}(w)$ we obtain
\begin{eqnarray*}
\mathbf{E} [\varphi(u(t,x))]&=&
 \int_{\mathbb{R}^{d}}  \ Ju_{0} \   \rho(u_{0}^{-1}(w),t,x)  \ \varphi(Z_t(w) )   dw \\
&=&
 \int_{\mathbb{R}^{d}} \  Ju_{0}( Z_t^{-1}(y)) \  JZ_{t}  \ \rho(u_{0}^{-1}( Z_t^{-1}(y) ),t,x) \  \varphi(y)  \  dy
\end{eqnarray*}
and thus relation (\ref{9a-2}) is obtained. \qed

\begin{remark}
The assumption $\div b=0$ can be interpreted as follows (see \cite{L}): in fluid mechanics or more generally in continuum mechanics, incompressible
flow (isochoric flow) refers to a flow in which the material density is
constant within a fluid parcel—an infinitesimal volume that moves with the
velocity of the fluid. This is equivalent to the condition that  the divergence of the fluid velocity is zero.

\end{remark}

\section{Appendix}

We present here some basic element on the fractional Brownian motion and on the Malliavin calculus.

\subsection{Fractional Brownian motion}

Consider $(B_{t}^{H})_{t\in\lbrack0,T]}$ a fractional Brownian
motion with Hurst  parameter $H\in(0,1)$. Recall that is it a centered Gaussian  process with covariance function 
\begin{equation}
\label{cov}
\mathbf{E} B ^{H}_{t} B ^{H}_{s} :=R_{H}(t,s)= \frac{1}{2} ( t^{2H} + s^{2H}- \vert t-s\vert ^{2H}, \hskip0.5cm s,t \in [0,T].
\end{equation}
The fractional Brownian motion can be also defined as the only self-similar Gaussian process with stationary increments.

Denote by ${\mathcal{H}}$ its
canonical Hilbert space . If $H=\frac{1}{2}$ then $B^{\frac{1}{2}}$ is the
standard Brownian motion (Wiener process) $W$ and in this case ${\mathcal{H}%
}=L^{2}([0,T])$. Otherwise $\mathcal{H}$ is the  Hilbert space  on $[0,T]$ extending the set of indicator function $\mathbf{1}_{[0,T]}, t\in [0,T]$  (by linearity and closure under the inner product) the rule
\[
\left\langle \mathbf{1}_{[0,s]};\mathbf{1}_{[0,t]}\right\rangle _{\mathcal{H}%
}=R_{H}\left(  s,t\right)  :=2^{-1}\left(  s^{2H}+t^{2H}-\left\vert
t-s\right\vert ^{2H}\right)  .
\]

 The followings facts will be needed in the sequel (we refer to
\cite{PiTa1} or \cite{N} for their proofs):
\begin{description}
\item{$\bullet$ } If $H>\frac{1}{2}$, the elements of ${\cal{H}}$ may be not functions
but distributions; it is therefore more practical to work with
subspaces of ${\cal{H}}$ that are sets of functions. Such a subspace
is
\begin{eqnarray*}  \left| {\cal{H}}\right| &=&\left \{ f:[0,T]\to
\mathbb{R} \,\,\Big |  \int _{0}^{T} \int_{0}^{T} \vert f(u)\vert
\vert f(v)\vert \vert u-v\vert ^{2H-2} dvdu <\infty \right \}.
\end{eqnarray*}
Then $\left| {\cal{H}}\right|$ is a strict subspace of $
{\cal{H}}$ and we actually have the inclusions
\begin{eqnarray}
\label{inclu1} L^{2}([0,T]) \subset
L^{\frac{1}{H} } ([0,T]) \subset \left| {\cal{H}}\right|
\subset  {\cal{H}}.
\end{eqnarray}
\item{$\bullet$ } The space $\left| {\cal{H}}\right|$ is not complete with respect to the norm $\Vert \cdot
\Vert _{{\cal{H}}}$ but it is a Banach  space with respect to the
norm
\begin{eqnarray*}
 \Vert f\Vert ^{2}_{\left| {\cal{H}}\right|
}&=&\int _{0}^{T} \int_{0}^{T} \vert f(u)\vert \vert
f(v)\vert \vert u-v\vert ^{2H-2} dvdu .
\end{eqnarray*}

\item{$\bullet$ } If $H>\frac{1}{2}$ and $f,g$ are two elements in the space  $\left| {\cal{H}}\right|$, their scalar product in ${\cal{H}}$ can be expressed by
    \begin{equation}
    \label{27i-1}
    \langle f,g\rangle _{{\cal{H}}}=\alpha _{H} \int_{0} ^{T} \int_{0} ^{T} dudv \vert u-v\vert ^{2H-2} f(u) g(v)
    \end{equation}
where $\alpha _{H}= H(2H-1)$.

\item{$\bullet$ } when $H<\frac{1}{2}$ then the canonical Hilbert space is a space of functions. We have
$$ C^{\gamma } \subset \mathcal{H} \subset L ^{2} ([0,T])$$
for all $\gamma > \frac{1}{2}-H$ where $C^{\gamma}$ denotes the class of H\"older continuous functions of order $\gamma$.

\item{$\bullet$ } The fBm \index{fractional Brownian motion} admits a representation as Wiener integral \index{Wiener integral}of the form
\begin{equation}
B ^{H}_{t}=\int_{0}^{t}K_{H}(t,s)dW_{s}, \label{BH}%
\end{equation}
where $W=\{W_{t},t\in T\}$ is a Wiener process, and $K_{ H}(t,s)$ is the kernel
\begin{equation}
K_{H}( t,s)=d_{H}\left(  t-s\right)  ^{H-\frac{1}{2}}+s^{H-\frac
{1}{2}}F_{1}\left(  \frac{t}{s}\right)  , \label{for1}%
\end{equation}
$d_{H}$ being a constant and
\begin{equation*}
F_{1}\left(  z\right)  =d_{H}\left(  \frac{1}{2}-H\right)  \int _{0}^{z-1}\theta^{H-\frac{3}{2}}\left( 1-\left(
\theta+1\right)
^{H-\frac{1}{2}}\right)  d\theta.
\end{equation*}
If $H>\frac{1}{2}$, the kernel $K_{H}$ has the simpler expression
\begin{equation}
 \label{K}
K_{H}(t,s)= c_{H} s^{\frac{1}{2}-H} \int _{s}^{t} (u-s)^{H-\frac{3}{2}} u^{H-\frac{1}{2}}  du
 \end{equation}
 where $t>s$ and $c_{H} =\left( \frac{ H(H-1) }{\beta( 2-2H, H-\frac{1}{2}) } \right)
^{\frac{1}{2}}.$
\end{description}

A a $d$ dimensional fractional Brownian motion $B^{H}= (B^{H_{1}}, \ldots, B ^{H_{d}}$ with Hurst parameter $H=(H_{1}, \ldots , H_{d} ) \in (0,1) ^{d}$ is a centered Gaussian process in $\mathbb{R} ^{d} $ with independent components and the covariance of  the $i$th  component is given by 
$$R_{H_{i}}(t,s)=\mathbf{E} B^{H_{i}}_{t}B^{H_{i}} _{s} =\frac{1}{2} ( t^{2H_{i}}+ s^{2H_{i}}-\vert t-s\vert ^{2H_{i}})
$$
for every $1\leq i\leq d$. 

\subsection{The Malliavin derivative}

Here we describe the elements from the Malliavin calculus  that we  need in the paper.  We refer \cite{N} for a more complete exposition. Consider ${\mathcal{H}}$ a real separable Hilbert space and $(B (\varphi), \varphi\in{\mathcal{H}})$ an isonormal Gaussian process \index{Gaussian process} on a probability space $(\Omega, {\cal{A}}, P)$, which is a centered Gaussian family of random variables such that $\mathbf{E}\left( B(\varphi) B(\psi) \right)  = \langle\varphi, \psi\rangle_{{\mathcal{H}}}$.

We denote by $D$  the Malliavin  derivative operator that acts on smooth functions of the form $F=g(B(\varphi _{1}), \ldots , B(\varphi_{n}))$ ($g$ is a smooth function with compact support and $\varphi_{i} \in {{\cal{H}}}, i=1,...,n$)
\begin{equation*}
DF=\sum_{i=1}^{n}\frac{\partial g}{\partial x_{i}}(B(\varphi _{1}), \ldots , B(\varphi_{n}))\varphi_{i}.
\end{equation*}
It can be checked that the operator $D$ is closable from $\mathcal{S}$ (the space of smooth functionals as above) into $ L^{2}(\Omega; \mathcal{H})$ and it can be extended to the space $\mathbb{D} ^{1,p}$ which is the closure of $\mathcal{S}$ with respect to the norm
\begin{equation*}
\Vert F\Vert _{1,p} ^{p} = \mathbf{E} F ^{p} + \mathbf{E} \Vert DF\Vert _{\mathcal{H}} ^{p}. 
\end{equation*}
We denote by  $ \mathbb{D} ^{k, \infty}:= \cap _{p\geq } \mathbb{D} ^{k,p}$ for every $k\geq 1$. In our paper, $\mathcal{H}$ will be the canonical Hilbert space associated with the fractional Brownian motion, as defined in the previous paragraph.


\begin{thebibliography}{99}


\bibitem{AB}
{O. Aboura and S.  Bourguin (2013): }{\em  Density estimates for solutions to one dimensional backward SDE's.  }vPotential Anal. 38(2), 573-587.

\bibitem{AMR}
{ R. Abraham, J. Marsden and T. Ratiu (1988): } {\em Manifolds, tensor analysis, and
applications. Second edition. } Applied Mathematical Sciences, 75.
Springer-Verlag, New York.

\bibitem{BKT}
{M. Besal\'u, A. Kohatsu-Higa and S. Tindel (2013): }{\em Gaussian type lower bounds for the density of solutions of SDEs driven by fractional Brownian motions. } Preprint.


\bibitem{CO}
{P. Catuogno and C. Oliveira (2013): }{\em $L ^{p}$ solutions of the stochastic transport equation. } Random Operators and Stochastic Equations,  21(2), 125-134. 

\bibitem{Chow}
 {P. L. Chow (2007): }     Stochastic Partial Differential Equations, Chapman  Hall/CRC, 2007.


\bibitem{DL}
{R. DiPerna,  P. L. Lions (1989): } {\em Ordinary differential
equations, transport theory and Sobolev spaces }, Invent. Math.,  98, 
511--547.


\bibitem{Fre1}
 F Fedrizzi , F. Flandoli. {\it Noise prevents singularities in linear transport equations},  Journal of Functional Analysis, 264,  1329-1354, 2013.





\bibitem{FGP2}
 F. Flandoli, M. Gubinelli, E. Priola, 2010. {\it Well-posedness of the transport equation by stochastic
 perturbation}, Invent. Math., 180(1):  1-53


\bibitem{FlandRusso}
{ F. Flandoli, F. Russo (2002): } {\em   Generalized integration and stochastic ODEs. }  Ann. Probab.  30(1), 270-292.


\bibitem{Ku}
{H. Kunita (1984): }{\em  Stochastic differential equations and stochastic flows of diffeomorphisms. } Ecole d'\'et\'e de probabilit\'es de Saint-Flour, XII—1982, 143–303, Lecture Notes in Math., 1097, Springer, Berlin.


\bibitem{Ku3}
{H. Kunita (1984): } {\em  First Order Stochastic Partial Differential Equations }, in Proceedings of
the Taniguchi International Symposium  on  Stochastic Analysis,  North-Holland Mathematical Library,  249-269.






\bibitem{Ku2}
 {H. Kunita (1990):}  Stochastic flows and stochastic differential
equations, Cambridge University Press, 1990.

\bibitem{L}
{P.L.Lions (1996): } {\em  Mathematical topics in fluid mechanics, Vol. I:
incompressible models. } Lecture Series in Mathematics and its
applications, 3, Oxford University Press.


\bibitem{Maurelli}
{M Maurelli (20011): } {\em  Wiener chaos and uniqueness for stochastic transport equation}, 
Comptes Rendus Mathematique,  349,  11-12,  669--672. 

\bibitem{MNP}
{S.E.A. Mohammed, T.  Nilssen and  F. Proske (2012): } {\em Sobolev Differentiable Stochastic Flows of SDE`s with Measurable Drift and Applications. } Preprint.


\bibitem{NV}
{I. Nourdin and F. Viens (2009): }{\em Density formula and concentration inequalities with Malliavin calculus. } Electronic Journal of Probability, 14, paper 78,  2287-2309.


\bibitem{N}
{D. Nualart (2006): }{\em Malliavin Calculus and Related Topics. Second Edition.  }{Springer New York. }

\bibitem{NQ}
{D. Nualart and L.  Quer-Sardanyons (2009): }{\em  Gaussian density estimates for solutions to quasi-linear stochastic partial differential equations. } Stochastic Process. Appl. 119(11), 3914-3938.

\bibitem{Oli}
  {C. Olivera (2014): } {\em  Well-posedness of   first order  semilinear  PDEs by stochastic
 perturbation. } Nonlinear Anal. 96, 211-215.



\bibitem{PiTa1}
{V. Pipiras and Murad Taqqu (2001): }{\em Integration questions
related to the fractional Brownian motion. } {Probability Theory and
Related Fields, 118, 2, 251-281. }


\bibitem{RV}
{F. Russo and P. Vallois (1993): }{\em Forward, backward and symmetric stochastic integration. } Probab. Theory Rel. Fileds, 97(3), 403-421. 

\bibitem{RV1}
{F. Russo and P. Vallois (2007): }{\em Elements of stochastic calculus via regularization. } S\'eminaire de Probabilit\'es XL, Lecture Notes in Mathematics 1899, 147-186.  

\bibitem{Sans}
{M. Sanz-Sol\'e (1995): } {\em Malliavin Calculus.  With applications to stochastic partial differential equations. } Fundamental Sciences, EPFL Press, Lausanne (2005). 


\end{thebibliography}
\end{document}